\input amstex
\documentstyle{amsppt}
\magnification 1200
\define\lra{\longrightarrow}
\define\p{\text{\rm Prob}}
\define\toD{\overset\Cal D\to\lra}
\define\toP{\overset p\to\lra}
\define\Z{\Bbb Z}
\define\E{\Bbb E}
\define\Y{\Bbb Y}
\define\la{\lambda}
\TagsOnRight
\NoRunningHeads
\topmatter
\address
Department of Mathematics, The University of
Pennsylvania, Philadelphia, PA 19104-6395, U.S.A.  E-mail address:
{\tt borodine\@math.upenn.edu}
\endaddress
\title
Longest increasing subsequences of random colored permutations
\endtitle
\abstract
We compute the limit distribution for (centered and scaled) length of the longest increasing subsequence of random colored permutations. The limit distribution function is a power of that for usual random permutations computed recently by Baik, Deift, and Johansson (math.CO/9810105). In two--colored case our method provides a different proof of a similar result by Tracy and Widom about longest increasing subsequences of signed permutations (math.CO/9811154).

Our main idea is to reduce the `colored' problem to the case of usual random permutations using certain combinatorial results and elementary probabilistic arguments. 

 \endabstract
\author
Alexei Borodin 
\endauthor
\endtopmatter
\head 1. Introduction
\endhead
Baik, Deift, and Johansson recently solved a problem about the asymptotic behavior of the length $l_n$ of the longest increasing subsequence for the random permutations of order $n$ as $n\to\infty$ (with the uniform distribution on the symmetric group $S_n$). They proved, see [BDJ], that the sequence $$\left\{\frac{l_n-2\sqrt{n}}{n^{1/6}}\right\}$$ converges in distribution, as $n\to\infty$, to a certain random variable whose distribution function we shall denote by $F(x)$. This distribution function can be expressed via a solution of the Painlev\'e II equation, see [BDJ] for details. It was first obtained by Tracy and Widom [TW1] in the framework of Random Matrix Theory  where it gives the limit distribution for the (centered and scaled) largest eigenvalue in the Gaussian Unitary Ensemble of Hermitian matrices.

The problem of the asymptotics of $l_n$ was first raised by Ulam [U]. Substantial contributions to the solution of the problem have been made by Hammersley [H], Logan and Shepp [LS], Vershik and Kerov [VK1, VK2].

 A survey of the interesting history  of this problem, further references, and a discussion of its intriguing connection with Random Matrix Theory can be found in [BDJ].

Soon after the appearance of [BDJ] Tracy and Widom computed the asymptotic behavior of the length $l_n'$ of the longest increasing for the random `signed permutations', see definitions in the next section. In [TW2] they showed that $$\left\{\frac{l'_n-2\sqrt{2n}}{2^{2/3}(2n)^{1/6}}\right\}$$ converges in distribution, as $n\to\infty$, to a random variable with the distribution function $F^2(x)$.

The present paper provides another proof of the result by Tracy and Widom. In our approach the distribution function $F^2(x)$ arises as the distribution function of the maximum of two asymptotically independent variables each of which behaves as $(l_n-2\sqrt{n})/n^{1/6}$ (hence, by [BDJ], converges to the distribution given by $F(x)$). 

The combinatorial techniques we use relies on recent works by Rains [R] and Fomin \& Stanton [FS]. It also allows to handle a more general case of `colored permutations' (the problem for `two--colored case', essentially, coincides with that for signed permutations). We show that for the length $l''_n$ of the longest increasing subsequence of the random $m$--colored permutations of order $n$ the sequence $$\left\{\frac{l''_n-2\sqrt{mn}}{m^{2/3}(mn)^{1/6}}\right\}$$ converges in distribution, as $n\to\infty$, to a random variable with distribution function $F^m(x)$. The function $F^m(x)$ naturally appears as the distribution function of the maximum of $m$ asymptotically independent variables, each having $F(x)$ as the limit distribution function.

Combinatorial quantities which we consider can be also interpreted as expectations of certain central functions on unitary groups, see Section 4.  

I am very grateful to G.~I.~Olshanski for a number of valuable discussions.

\head
2. Colored permutations and signed permutations
\endhead
A {\it colored permutation} is a map from $\{1,\dots,n\}$ to $\{1,\dots,n\}\times\{1,\dots,m\}$ such that its composition with the projection on the first component of the target set is a permutation (of order $n$). One can view such a map as a permutation with one of $m$ colors attributed to each of $n$ points which this permutation permutes. The set of all colored permutations of order $n$ with $m$ colors will be denoted by $S_n^{(m)}$.

An {\it increasing subsequence} of $\pi\in S_n^{(m)}$ is a sequence $1\le i_1<i_2<\dots<i_k\le n$ such that the first coordinates of $\pi(i_j)$ increase in $j$ and the second coordinates of $\pi(i_j)$ are equal. Thus, elements of an increasing subsequence are of the same color, say, $p$.
The {\it length} of such increasing subsequence is defined to be $m(k-1)+p$.

These definitions are due to Rains [R]. A slightly more general notion of {\it hook permutation} was introduced and intensively used earlier by Stanton and White [SW].

We shall consider $S_n^{(m)}$ as a probability space with uniform distribution: probability of each colored permutation  is $|S_n^{(m)}|^{-1}=(m^{n}n!)^{-1}$. Then the length of the longest increasing subsequence becomes a random variable on this space, it will be denoted as $L_n^{{col}(m)}$.

Let $H_n$ be the hyperoctahedral group of order $n$ defined as the wreath product $\Z_2^n\ltimes S_n$ ($S_n$ is the symmetric group of order $n$). The elements of $H_n$ are called {\it signed permutations}. This group can be naturally embedded in $S_{2n}$ as the group of permutations $\sigma$ of $\{-n,-n+1,\dots,-1,1,\dots,n-1,n\}$ subject to the condition
$\sigma(-x)=-\sigma(x)$. Indeed, each such permutation is parametrized by the permutation $|\sigma|\in S_n$ and the set of signs of $\sigma(1),\dots,\sigma(n)$. Using the natural ordering on the set $\{-n,-n+1,\dots,-1,1,\dots,n-1,n\}$ we can define the length of the longest increasing subsequence for each signed permutation. Assuming that every signed permutation has probability $|H_n|^{-1}=(2^{n}n!)^{-1}$, we get a random variable on $H_n$ which will be denoted as $L_n^{even}$.

The group $H_n$ can also be embedded into the symmetric group of order $2n+1$: we add 0 to the set $\{-n,-n+1,\dots,-1,1,\dots,n-1,n\}$ and assume that the elements $\sigma\in H_n$ satisfy the same condition $\sigma(-x)=-\sigma(x)$. Clearly, this implies $\sigma(0)=0$.
The random variable on $H_n$ equal to the length of the longest increasing subsequence with respect to this realization will be denoted by $L_n^{odd}$.

Note that for any element $\sigma\in H_n$, 
$$
L_n^{odd}(\sigma)-L_n^{even}(\sigma)=0\text{ or }1.
\tag 2.1
$$

\head 3. Rim hook tableaux \endhead

We refer to the work [SW] for the definitions concerning rim hook tableaux.

The next claim is a direct consequence of the Schensted algorithm, see [S].

\proclaim{Proposition 3.1} Permutations of order $n$ with the length of longest increasing subsequence equal to $l$ are in one--to--one correspondence with pairs of standard Young tableaux of the same shape with $n$ boxes and width $l$. 
\endproclaim

Here is a generalization of this claim for colored permutations. 

\proclaim{Proposition 3.2 [R], [SW]} Colored permutations with $m$ colors of order $n$ with the length of longest increasing subsequence equal to $l$ are in one--to--one correspondence with pairs of $m$--rim hook tableaux of the same shape with $mn$ boxes and width $l$.
\endproclaim

In [SW] it was proved that $\lceil\frac lm\rceil=\lceil \frac wm\rceil$ where $w$ is the width of the rim hook tableau corresponding to a permutation with the length of longest increasing subsequence equal to $l$ ($\lceil a\rceil$ stands for the smallest integer $\ge a$). The refinement of this statement given above was published in [R].

\proclaim{Proposition 3.3 [R]} Signed permutations of order $n$ embedded in $S_{2n}$ with the length of longest increasing subsequence equal to $l$ are in one--to--one correspondence with pairs of $2$--rim hook tableaux of the same shape with $2n$ boxes and width $l$.
\endproclaim

The length of the longest increasing subsequence for signed permutations embedded into the symmetric group of odd order can also be interpreted in terms of rim hook tableaux, see [R, proof of Theorem 2.3].

Note that Propositions 3.2 and 3.3 imply that the distributions of random variables $L^{col(2)}_n$ and $L^{even}_n$ coincide. 

\head 4. Expectations over unitary groups \endhead

Everywhere below the symbol $\E_{U\in U(k)}f(U)$ stands for the integral of $f$ over $U\in U(k)$ with respect to the Haar measure on the unitary group $U(k)$ normalized so that $\E_{U\in U(k)}1=1$ (i.e., $\E$ denotes the expectation of $f$ with respect to the uniform distribution on the unitary group).

\proclaim{Proposition 4.1 [R]} 
$$
\p\{L^{col(m)}_n\le k\}=(m^nn!)^{-1}\cdot\E_{U\in U(k)}\left(|{Tr(U^m)}^n|^2\right).
\tag 4.1
$$
\endproclaim

\proclaim{Proposition 4.2 [R]} 
$$
\p\{L^{even}_n\le k\}=(2^nn!)^{-1}\cdot\E_{U\in U(k)}\left(|{Tr(U^2)}^n|^2\right).
\tag 4.2
$$
$$
\p\{L^{odd}_n\le k\}=(2^nn!)^{-1}\cdot\E_{U\in U(k)}\left(|{Tr(U^2)}^nTr(U)|^2\right).
\tag 4.3
$$
\endproclaim
[DS] gives (4.1) for $k\ge mn$, (4.2) for $k\ge 2n$, and (4.3) for $k\ge 2n+1$. For such values of $k$ the left--hand sides of (4.1), (4,2), (4.3) are all equal to 1.

\head 5. Rim hook lattices \endhead
Our main reference for this section is the work [FS] by Fomin and Stanton. 

For this section we fix an integer number $m$, all our rim hooks here will contain exactly $m$ boxes.

Let $\mu$ and $\lambda$ be shapes (Young diagrams) such that
$\mu\subset\la$ and $\la-\mu$ is a ($m$--)rim hook. Then we shall write $\mu\nearrow\la$. 

We introduce a partial order on the set of Young diagrams as follows: $\la\succeq\mu$ if and only if there exists a sequence $\nu_1,\nu_2,\dots,\nu_k$ of Young diagrams such that
$\mu\nearrow\nu_1\nearrow\nu_2\nearrow\dots\nearrow\nu_k\nearrow\la$.
The empty Young diagram is denoted by $\emptyset$. We shall say that a Young diagram $\la$ is $m$--decomposable if $\la\succeq\emptyset$.

The poset of all $m$--decomposable shapes with $\succeq$ as the order is called {\it rim hook lattice} and is denoted by $RH_m$. (It can be shown that this poset is indeed a lattice).

For $m=1$ we get the {\it Young lattice}: the poset of all Young diagrams ordered by inclusion. The Young lattice will be denoted by $\Y$.

\proclaim{Proposition 5.1 [FS]} The rim hook lattice $RH_m$ is isomorphic to the Cartesian product of $m$ copies of the Young lattice: $RH_m\cong \Y^m$.
\endproclaim

In other words, $RH_m$ is isomorphic to the poset of $m$--tuples of Young diagrams with the following coordinate--wise ordering: 
one tuple is greater than or equal to another tuple if the $k$th  coordinate of the first tuple includes (i.e., greater than or equal to) the $k$th coordinate of the second tuple for all $k=1,\dots,m$.

Clearly, the number of $m$--rim hook tableaux of a given shape $\la$ is equal to the number of paths $\emptyset\nearrow\nu_1\nearrow\nu_2\nearrow\dots\nearrow\nu_k\nearrow\la$, $k={|\la|}/{m}-1$, from $\emptyset$ to $\la$ (and is equal to 0 if $\la$ is not $m$--decomposable), $|\la|$ stands for the number of boxes in $\la$. We shall denote this number by $\dim_m\la$ and call it the {\it $m$--dimension} of the shape $\la$.

Take any $\la\in RH_m$ and the corresponding $m$--tuple $(\la_1,\dots,\la_m)\in \Y^m$. Note that $|\la|=m(|\la_1|+\dots+|\la_m|)$. We have
$$
\dim_m\la=\frac {(|\la_1|+\dots+|\la_m|)!}{|\la_1|!\cdots|\la_m|!}\cdot \dim_1\la_1\cdots \dim_1\la_m.
\tag 5.1
$$

Indeed, to specify the path from $\emptyset$ to $(\la_1,\dots,\la_m)$ we need to specify $m$ paths from $\emptyset$ to $\la_k$ in the $k$th copy of $\Y$ for $k=1,\dots,m$ together with the order in which we make steps along those paths. The number of different orders is the combinatorial coefficient in the right--hand side of (5.1) while the number of different possibilities for the $m$ paths in $\Y$ is the product of 1--dimensions $\dim_1\la_1\cdots \dim_1\la_m$. 

Note that the number of pairs of $m$--rim hook tableaux of the same shape $\la$ is exactly 
$$
\dim^2_m\la=\left(\frac {(|\la_1|+\dots+|\la_m|)!}{|\la_1|!\cdots|\la_m|!}\right)^2\cdot \dim^2_1\la_1\cdots \dim^2_1\la_m.
\tag 5.2
$$

Let us denote by $w(\la)$ the width of a Young diagram $\la$. We shall need the following
\proclaim{Observation 5.2} For any $\la\in RH_m$ and corresponding $m$--tuple $(\la_1,\dots,\la_m)\in \Y^m$ we have
$$
m\cdot\max\{w(\la_1),\dots, w(\la_m)\}-w(\la)\in\{0,1,\dots,m-1\}.
\tag 5.3
$$
\endproclaim
 This follows immediately from the explicit construction of the isomorphism from Proposition 5.1, see [FS, \S2].
 
 Observation 5.2 will be crucial for our further considerations.

\head 6. Plancherel distributions \endhead

Using the correspondence from Proposition 3.2 (which is exactly the rim hook generalization of the Schensted algorithm, see [SW]), we can associate to each $m$--colored permutation 
of order $n$ a Young diagram with $mn$ boxes --- the common shape of the corresponding pair of $m$--rim hook tableaux.
The image of the uniform distribution on $S_n^{(m)}$ under this map gives a probability distribution on $m$--decomposable Young diagrams with $mn$ boxes; the weight of a Young diagram $\la$ is, clearly, equal to
$(m^nn!)^{-1}\cdot \dim_m^2\la$. As a consequence, we get
(cf. [FS, Corollary 1.6]) 
$$
\sum_{|\la|=mn} \dim_m^2\la=m^nn!.
\tag 6.1
$$

Using the isomorphism of Proposition 5.1 we transfer our probability distribution to the set of $m$--tuples of Young disgrams with total number of boxes equal to $n$. Then by (5.2) we see that the probability of an $m$--tuple 
$(\la_1,\dots,\la_m)\in \Y^m$ with $|\la_1|+\dots+|\la_m|=n$ equals
$$
\p\{(\la_1,\dots,\la_m)\}=\frac 1{m^nn!}\left(\frac {n!}{|\la_1|!\cdots|\la_m|!}\right)^2\cdot \dim^2_1\la_1\cdots \dim^2_1\la_m.
\tag 6.2
$$
This distribution will be called the {\it Plancherel distribution}.

\proclaim{Proposition 6.1} For any $n=1,2,\ldots$ 
$$
\gathered
\p\{|\la_1|=n_1,\dots|\la_m|=n_m;n_1+\dots+n_m=n\}\\
=\frac 1{m^n}\,\frac {n!}{n_1!\cdots n_m!}.
\endgathered
\tag 6.3
$$
\endproclaim
\demo{Proof} Direct computation. Using (6.1) for $m=1$ and (6.2) we get
$$
\gathered
\p\{|\la_1|=n_1,\dots|\la_m|=n_m;n_1+\dots+n_m=n\}\\
=\frac 1{m^nn!}\sum_{\Sb |\la_k|=n_k\\k=1,\dots,m\endSb}\left(\frac {n!}{|\la_1|!\cdots|\la_m|!}\right)^2\cdot \dim^2_1\la_1\cdots \dim^2_1\la_m\\
=\frac 1{m^nn!}\left(\frac {n!}{n_1!\cdots n_m!}\right)^2\cdot n_1!\cdots n_m!=\frac 1{m^n}\,\frac {n!}{n_1!\cdots n_m!}.\qed
\endgathered
$$
\enddemo

\head 7. Two lemmas from Probability Theory \endhead

We shall denote by $\toP$ convergence of random variables in probability, and by $\toD$ convergence in distribution. 

\proclaim{Lemma 7.1 [B, Theorem 4.1]}
Let random variables $\xi$, $\{\xi_n\}_{i=1}^\infty$, $\{\eta_n\}_{i=1}^\infty$ satisfy 
$$
\xi_n-\eta_n\toP 0,\quad \xi_n\toD \xi.
$$
Then
$$
\eta_n\overset\Cal D\to\lra \xi.
$$
\endproclaim

For $m$ real random variables $\xi_1,\dots,\xi_m$ we shall denote by $\xi_1\times\xi_2\times\cdots\times\xi_m$ a $\Bbb R^m$--valued random variable with distribution function
$$
F_{\xi_1\times\cdots\times\xi_m}(x)=F_{\xi_1}(x_1)\cdots F_{\xi_m}(x_m).
$$
\proclaim{Lemma 7.2}
Let $\{\xi^{(k)}_n\}_{n=1}^\infty$, $k=1,\dots,m$,  be $m>1$ sequences of real random variables convergent in distribution to random variables $\xi^{(1)},\dots,\xi^{(m)}$, respectively. Denote by $B^{(m)}_n$ the $n$th order $m$--dimensional fair Bernoulli distribution:
$$
\p\{B^{(m)}_n=(k_1,\dots,k_{m})\}=\cases\frac 1{m^n}\,\frac {n!}{k_1!\cdots k_m!},& k_1+\ldots+k_m=n,\ k_i\in\{0,\dots,n\}\\
0,&\text{otherwise}
\endcases
$$
Then the sequence 
$$
\zeta_n=\xi^{(1)}_{B^{(m)}_{n,1}}\times\xi^{(2)}_{B^{(m)}_{n,2}}\times\cdots\times\xi^{(m)}_{B^{(m)}_{n,m}}
$$
converges in distribution to $\xi^{(1)}\times\xi^{(2)}\times\dots\times\xi^{(m)}$. 
In particular, $\xi^{(1)}_{B^{(m)}_{n,1}}\toD \xi^{(1)}$.
\endproclaim
\demo{Proof}
The conditions $\xi^{(k)}_n\toD \xi^{(k)}$, $k=1,\dots,m$, are equivalent to the pointwise convergence of distribution functions
$$
F_{\xi^{(k)}_n}(x)\to F_{\xi^{(k)}}(x).
$$
We have
$$
F_{\zeta_n}(x)=\sum_{\Sb k_1+\dots+k_m=n\\  k_i=1,\dots,n\endSb}\frac 1{m^n}\,\frac{n!}{k_1!\cdots k_m!} \cdot F_{\xi^{(1)}_{k_1}}(x_1)\cdots F_{\xi^{(m)}_{k_m}}(x_m).
$$

For $m$ convergent number sequences $\{a^{(k)}_n\}$ with limits $a^{(k)}$, $k=1,\dots,m$,  the sequence 
$$
c_n=\sum_{\Sb k_1+\dots+k_m=n\\  k_i=1,\dots,n\endSb}\frac 1{m^n}\,\frac{n!}{k_1!\cdots k_m!}\,a^{(1)}_{k_1}\cdots a^{(m)}_{k_m}
$$ converges to the product $a^{(1)}\cdots a^{(m)}$, so we conclude, that $\{F_{\zeta_n}(x)\}$ converges to $F_\xi^{(1)}(x_1)\cdots F_\xi^{(m)}(x_m)$ pointwise. \qed
\enddemo

Lemma 7.2 implies that random variables $\xi^{(1)}_{B^{(m)}_{n,1}},\ \xi^{(2)}_{B^{(m)}_{n,2}},\dots,\ \xi^{(m)}_{B^{(m)}_{n,m}}$ are {\it asymptotically independent} as $n\to\infty$.

\head 8. Asymptotics \endhead

Baik, Deift, and Johansson recently proved, see [BDJ], that the sequence 
$$
\frac{L^{col(1)}_n-2\sqrt{n}}{n^{1/6}}
$$
converges in distribution to a certain random variable, whose distribution function will be denoted by $F(x)$. This function can be expressed through a particular solution of the Painlev\'e II equation, see [BDJ] for details.

In this section we shall study the asymptotic behavior of $L^{col(m)}_n$ for $m\ge 2$ when $n\to\infty$. Our main result is the following statement.
\proclaim{Theorem 8.1} For any $m=2,3,\dots$ the sequence
$$
\frac{L^{col(m)}_n-2\sqrt{mn}}{(mn)^{1/6}}
$$
converges in distribution, as $n\to\infty$, to a random variable with distribution function $F^m(m^{-\frac 23}x)$.
\endproclaim

This result for $m=2$ was proved by Tracy and Widom in [TW2]. 
Since the distributions of $L_n^{col(2)}$ and $L_n^{even}$ coincide (see Section 3), we immediately get two other asymptotic formulas.

\proclaim{Corollary 8.2 [TW2]} The sequence
$$
\frac{L^{even}_n-2\sqrt{2n}}{(2n)^{1/6}}
$$
converges in distribution, as $n\to\infty$, to a random variable with distribution function $F^2(2^{-\frac 23}x)$.
\endproclaim
\proclaim{Corollary 8.3 [TW2]} The sequence
$$
\frac{L^{odd}_n-2\sqrt{2n}}{(2n)^{1/6}}
$$
converges in distribution, as $n\to\infty$, to a random variable with distribution function $F^2(2^{-\frac 23}x)$.
\endproclaim
\demo{Proof of Corollary 8.3}
Relation (2.1) implies that with probability 1,
$$
(L^{odd}_n-L^{even}_n)\,(2n)^{-1/6}\to 0.
$$
Since the convergence with probability 1 implies the convergence in probability, the claim follows from Lemma 7.1 and Corollary 8.2.\qed
\enddemo
\demo{Proof of Theorem 8.1} The proof will consist of 4 steps.

\noindent{\bf Step 1.} Following Sections 3, 5, 6, we shall interpret $L^{col(m)}_n$ as a random variable on the space of all $m$--tuples $(\la_1,\dots,\la_m)$ of Young diagram with total number of boxes equal to $n$ supplied with the Plancherel distribution. Let us denote by $l_1^{(n)},\dots,l_m^{(n)}$ the widths of the Young diagrams $\la_1,\dots,\la_m$: 
$$
l_k^{(n)}=w(\la_k),\quad k=1,\dots,m.
$$

Observation 5.2 implies that with probability 1
$$
\frac{L^{col(m)}_n-m\cdot\max\{l_1^{(n)},\dots,l_m^{(n)}\}}{n^{1/6}}\to 0.
$$
Hence, by Lemma 7.1, it is enough to prove that 
$$
\p\left\{\frac{m\cdot\max\{l_1^{(n)},\dots,l_m^{(n)}\}-2\sqrt{mn}}{(mn)^{1/6}}\le x\right\}\to F^m(m^{-\frac 23}x)
$$

Our strategy is to prove that $l_1^{(n)},\dots,l_m^{(n)}$ are asymptotically independent, and that each of them asymptotically behaves as $L_{[n/m]}^{col(1)}$. Then Theorem 8.1 will follow from the result of [BDJ] stated  in the beginning of this section.

\noindent{\bf Step 2.}
$$
\gathered
\p\left\{\frac{m\cdot\max\{l_1^{(n)},\dots,l_m^{(n)}\}-2\sqrt{mn}}{(mn)^{1/6}}\le x\right\}\\
=\p\left\{\frac{m\cdot l_1^{(n)}-2\sqrt{mn}}{(mn)^{1/6}}\le x;\dots;\frac{m\cdot l_m^{(n)}-2\sqrt{mn}}{(mn)^{1/6}}\le x\right\}\\
=\p\left\{\frac{l_1^{(n)}-2\sqrt{n/m}}{(n/m)^{1/6}}\le m^{\frac 23}x;\dots;\frac{ l_m^{(n)}-2\sqrt{n/m}}{(n/m)^{1/6}}\le m^{\frac 23}x\right\}.
\endgathered
$$
Thus, it suffices to prove that
$$
\p\left\{\frac{l_1^{(n)}-2\sqrt{n/m}}{(n/m)^{1/6}}\le x_1;\dots;\frac{ l_m^{(n)}-2\sqrt{n/m}}{(n/m)^{1/6}}\le x_m\right\}\to F(x_1)\cdots F(x_m).
\tag 8.1
$$

\noindent{\bf Step 3.} Denote by $n_i^{(n)}$ the number of boxes in $\la_i$, $i=1,\dots,m$. Then $n_1^{(n)}+\dots+n_m^{(n)}=n$. We claim that for all $i=1,\dots,m$
$$
\frac{l_i^{(n)}-2\sqrt{n/m}}{(n/m)^{ 1/6}}-\frac{l_i^{(n)}-2\sqrt{n_i^{(n)}}}{{(n_i^{(n)})}^{1/6}}\toP 0.
\tag 8.2
$$ 
Proposition 6.1 implies that 
$$
\p\{n_i^{(n)}=k\}=\frac {(m-1)^{n-k}}{m^n}\,\binom{n}{k}.
$$
This means that $n_i^{(n)}$ can be interpreted as the sum 
$\xi_1+\dots+\xi_n$ of $n$ independent identically distributed Bernoulli variables with 
$$
\p\{\xi_j=0\}=\frac{m-1}m,\quad \p\{\xi_j=1\}=\frac{1}m.
$$
The central limit theorem implies that with probability converging to 1,  
$$
\frac nm-\left(\frac nm\right)^{\frac 12 +\varepsilon}<n_i^{(n)}< \frac nm+\left(\frac nm\right)^{\frac 12 +\varepsilon}
$$ 
for any $\varepsilon>0$.
Hence,
$$
\gathered
\left|{\left(\frac nm\right)}^{\frac 13}-({n_i^{(n)}})^{\frac 13}\right|< \left|\left(\frac nm\right)^{\frac 13}-\left(\frac nm\pm \left(\frac nm\right)^{\frac 12 +\varepsilon}\right)^{\frac 13}\right|\\=\left(\frac nm\right)^{\frac 13}\left|1-\left(1\pm \left(\frac nm\right)^{-\frac 12 +\varepsilon}\right)^{\frac 13}\right|<const\cdot \left(\frac nm\right)^{-\frac 16+\varepsilon}
\endgathered
\tag 8.3
$$
for sufficiently large $n$.

Similarly, we have
$$
\gathered
\left|{\left(\frac nm\right)}^{-\frac 16}-({n_i^{(n)}})^{-\frac 16}\right|< \left|{\left(\frac nm\right)}^{-\frac 16}-\left(\left(\frac nm\right)\pm \left(\frac nm\right)^{\frac 12 +\varepsilon}\right)^{-\frac 16}\right|\\={\left(\frac nm\right)}^{-\frac 16}\left|1-\left(1\pm \left(\frac nm\right)^{-\frac 12 +\varepsilon}\right)^{-\frac 16}\right|<const\cdot \left(\frac nm\right)^{-\frac 12-\frac 16+\varepsilon}<const\cdot \left(\frac nm\right)^{-\frac 12-\varepsilon}
\endgathered
\tag 8.4
$$
if we choose $\varepsilon<1/12$.

The result of [BDJ] implies that $L_n^{col(1)}n^{-\frac 12-\delta}\toD 0$ for any $\delta>0$.

From Proposition 6.2 it follows that in the notation of Lemma 7.2
$$
\left\{l_i^{(n)}\right\}_{n=1}^\infty=\left\{L_{(B^{(m)}_{n,1})}^{col(1)}\right\}_{n=1}^\infty.
$$
Applying Lemma 7.2 we see that
$$
\frac {l_i^{(n)}}{\left(n_i^{(n)}\right)^{\frac 12 +\delta}}\toD 0,
$$
which means that $l_i^{(n)}<\left(n_i^{(n)}\right)^{\frac 12 +\delta}$ with probability converging to 1.

Since $n_i^{(n)}<n/m+\left(n/m\right)^{\frac {1+\varepsilon}2}$ with probability converging to 1, we get that 
$$
\gathered
l_i^{(n)}<\left(\frac nm+\left(\frac nm\right)^{\frac {1+\varepsilon}2}\right)^{\frac 12 +\delta}=
\left(\frac nm\right)^{\frac 12 +\delta}\left(1+\left(\frac nm\right)^{\frac {-1+\varepsilon}2}\right)^{\frac 12 +\delta}\\<\left(\frac nm\right)^{\frac 12+\delta}+const\cdot \left(\frac nm\right)^{\delta+\frac \varepsilon2}<const\cdot \left(\frac nm\right)^{\frac 12+\delta}
\endgathered
\tag 8.5
$$
with probability converging to 1.

Making use of the relations (8.3), (8.4), (8.5), we now see that
$$
\gathered
\left|\frac{l_i^{(n)}-2\sqrt{n/m}}{(n/m)^{ 1/6}}-\frac{l_i^{(n)}-2\sqrt{n_i^{(n)}}}{{(n_i^{(n)})}^{1/6}}\right|\le l_i^{(n)}\left|{\left(\frac nm\right)}^{-\frac 16}-({n_i^{(n)}})^{-\frac 16}\right|\\+2\left|{\left(\frac nm\right)}^{\frac 13}-({n_i^{(n)}})^{\frac 13}\right|<const\cdot \left(\frac nm\right)^{\frac 12+\delta}\left(\frac nm\right)^{-\frac 12-\varepsilon}+const\cdot \left(\frac nm\right)^{-\frac 16+\varepsilon}
\endgathered
$$
with probability converging to 1. Since the last expression converges to zero as $n\to\infty$ if $\delta<\varepsilon<1/6$, the proof of (8.2) is complete.

The relation (8.2) implies that we have asymptotic equivalence of two $m$--di\-mensional vectors
$$
\gathered
\left(\frac{l_1^{(n)}-2\sqrt{n/m}}{(n/m)^{1/6}};\dots;\ \frac{l_m^{(n)}-2\sqrt{n/m}}{(n/m)^{1/6}}\right) \\-\left(\frac{l_1^{(n)}-2\sqrt{n_1^{(n)}}}{{(n_1^{(n)})}^{ 1/6}};\dots;\frac{l_m^{(n)}-2\sqrt{n_m^{(n)}}}{{(n_m^{(n)})}^{ 1/6}}\right)\toP 0.
\endgathered
\tag 8.6
$$

\noindent {\bf Step 4.} The random $m$--dimensional vector
$$
\zeta_{n}=\left(\frac{l_1^{(n)}-2\sqrt{n_1^{(n)}}}{{(n_1^{(n)})}^{1/6}};\dots;\ \frac{l_m^{(n)}-2\sqrt{n_m^{(n)}}}{{(n_m^{(n)})}^{ 1/6}}\right)
\tag 8.7
$$
is obtained from $m$ identical 1--dimensional random variables 
$$
\xi^{(i)}_n=\frac{L_n^{col(1)}-2\sqrt{n}}{{n}^{ 1/6}},\quad i=1,\dots,m 
$$
by the procedure described in Lemma 7.2:
$$
\zeta_{n}=\xi^{(1)}_{B^{(m)}_{n,1}}\times\xi^{(2)}_{B^{(m)}_{n,2}}\times\cdots\times\xi^{(m)}_{B^{(m)}_{n,m}}
$$
(this follows from Proposition 6.1).  Since sequences $\{\xi^{(i)}_n\}$ converge in distribution, as $n\to\infty$, to a random variable with distribution function $F(x)$ (this is the result of [BDJ]), the sequence (8.7) also converges in distribution to a random variable with the  distribution function $F(x_1)\cdots F(x_m)$. Then Lemma 7.1 and (8.6) conclude the proof of (8.1).\qed
\enddemo

\enddocument